\documentclass[12pt,twoside]{article} 
\usepackage[english]{babel}
\usepackage[latin1]{inputenc} 
\usepackage{amsmath}
\usepackage{amssymb,amsfonts}
\usepackage{graphicx}                   
\usepackage{times,amssymb,amscd} 

\newcommand{\bA}{\mathbf{A}}

\newcommand{\bE}{\mathbf{E}}
\newcommand{\bG}{\mathbf{G}}

\newcommand{\bR}{\mathbf{R}}
\newcommand{\bO}{\mathbf{O}}
\newcommand{\bS}{\mathbf{S}}
\newcommand{\bV}{\mathbf{V}}

\newcommand{\be}{\mathbf{e}}

\newcommand{\bx}{\mathbf{x}}
\newcommand{\by}{\mathbf{y}}

\newcommand{\bg}{\mathbf{g}}

\newcommand{\bI}{\mathbf{I}}

\newcommand{\BV}{\boldsymbol{V}}

\newcommand{\Be}{\boldsymbol{e}}
\newcommand{\Bu}{\boldsymbol{u}}
\newcommand{\Bv}{\boldsymbol{v}}

\newcommand{\cD}{\mathcal{D}}

\newcommand{\cP}{\mathcal{P}}
\newcommand{\cS}{\mathcal{S}}

\newcommand{\cB}{\mathcal{B}}

\newcommand{\SXR}{\bS^2\!\times\!\bR}

\newcommand{\NIL}{\mathbf{Nil}}
\newcommand{\SOL}{\mathbf{Sol}}

\begin{document}
\pagestyle{myheadings}
\markboth{\centerline{Jen\H o Szirmai}}
{Geodesic ball packings in $\SXR$ space}
\title
{Simply transitive geodesic ball packings to $\SXR$ space groups generated by glide reflections \footnote{AMS Classification 2010: 52C17, 52C22, 53A35, 51M20}}

\author{Jen\H o Szirmai \\
\normalsize Budapest University of Technology and \\ 
\normalsize Economics, Institute of Mathematics, \\
\normalsize Department of Geometry \\
\normalsize szirmai@math.bme.hu\\
\date{\normalsize{\today}}}


\maketitle
\begin{abstract}

The $\SXR$ geometry can be derived by the direct product of the spherical plane $\bS^2$ and the real line $\bR$.
In \cite{F01} J.~Z. {\sc Farkas} has classified and given the complete list of the space groups of $\SXR$.
The $\SXR$ manifolds were classified by E. {\sc Moln\'ar} and J.~Z. {\sc Farkas} in \cite{FM01} by similarity and diffeomorphism.
In \cite{Sz11} we have studied the geodesic balls and their volumes in $\SXR$ space, moreover we have introduced the notion of geodesic ball packing 
and its density and have determined the densest geodesic ball packing for generalized Coxeter space groups of $\SXR$.

In this paper we study the locally optimal ball packings to the $\SXR$ space groups having Coxeter point groups 
and at least one of the generators is a glide reflection. 
We determine the densest simply transitive geodesic ball arrangements for the above space groups, moreover we compute their optimal 
densities and radii.

The density of the densest packing is $\approx 0.80407553$, may be surprising enough
in comparison with the Euclidean result $\frac{\pi}{\sqrt{18}} \approx 0.74048$. E. {\sc Moln\'ar} has shown in \cite{M97}, that the homogeneous 3-spaces
have a unified interpretation in the real projective 3-sphere $\mathcal{PS}^3(\bV^4,\BV_4, \mathbb{R})$. 
In our work we shall use this projective model of $\SXR$ geometry.

\end{abstract}

\newtheorem{Theorem}{Theorem}[section]
\newtheorem{corollary}[Theorem]{Corollary}
\newtheorem{lemma}[Theorem]{Lemma}
\newtheorem{exmple}[Theorem]{Example}
\newtheorem{definition}[Theorem]{Definition}
\newtheorem{rmrk}[Theorem]{Remark}
\newtheorem{proposition}[Theorem]{Proposition}
\newenvironment{remark}{\begin{rmrk}\normalfont}{\end{rmrk}}
\newenvironment{example}{\begin{exmple}\normalfont}{\end{exmple}}
\newenvironment{acknowledgement}{Acknowledgement}



\section{Introduction}

$\SXR$ is derived as the direct product of the spherical plane $\bS^2$ and the real line $\bR$. 
The points are described by $(P,p)$ where $P\in \bS^2$ and $p\in \bR$. The isometry group $Isom(\SXR)$ of $\SXR$ can be derived by the direct product 
of the isometry group of the sphere $Isom(\bS^2)$ and the isometry group of the real line $Isom(\bR)$. 
\begin{equation}
\begin{gathered}
Isom(\bS^2):=\{A \in \bO(3) \ : \ \bS^2 \mapsto \bS^2 \ : \ (P,p) \mapsto (PA,p) \} \ \text{for any fixed $p$}.  \\
Isom(\bR):=\{\rho ~ : ~ (P,p) \mapsto (P, \pm p + r) \}, \ \text{for any fixed $P$}.\\
\text{here the "-" sign provides a reflection in the point} \ \frac{r}{2} \in \bR, \\ \text{by the "+" sign we get a translation of $\bR$}.  
\end{gathered} \tag{1.1}
\end{equation}
The structure of an isometry group $\Gamma \subset Isom(\SXR)$  is the following: $\Gamma:=\{(A_1 \times \rho_1), \dots (A_n \times \rho_n) \}$, where
$A_i \times \rho_i:=A_i \times (R_i,r_i):=(g_i,r_i)$, $(i \in \{ 1,2, \dots n \}$ and $A_i \in Isom(\bS^2)$, $R_i$ is either the identity map 
$\mathbf{1_R}$ of $\bR$ or the point reflection $\overline{\mathbf{1}}_{\mathbf{R}}$. $g_i:=A_i \times R_i$ is called the linear part of the transformation
$(A_i \times \rho_i)$ and $r_i$ is its translation part. 
The multiplication formula is the following:
\begin{equation}
(A_1 \times R_1,r_1) \circ (A_2 \times R_2,r_2)=((A_1A_2 \times R_1R_2,r_1R_2+r_2). \tag{1.2}
\end{equation}
\begin{definition}
A group of isometries $\Gamma \subset Isom(\SXR)$ is called {\it space group} if the linear parts form a finite group $\Gamma_0$ called the point group of 
$\Gamma$, moreover, the translation parts to the identity of this point group are required to form a one dimensional lattice $L_{\Gamma}$ of $\bR$.
\end{definition}
\begin{rmrk}
\begin{enumerate}
\item It can be proved that the space group $\Gamma$ has a compact fundamental domain $\mathcal{F}_\Gamma$.  
\item If $\Gamma$ is not assumed to have a lattice, then it may have an infinite point group $\Gamma_0$.
\end{enumerate}
\end{rmrk}
\begin{definition}
The $\SXR$ space groups $\Gamma_1$ and $\Gamma_2$ are geometrically equivalent, called equivariant, if there is a "similarity" transformation
$\Sigma:=S \times \sigma$ $(S \in Isom(\bS^2), \sigma \in Sim(\bR))$, such that $\Gamma_2=\Sigma^{-1} \Gamma_1 \Sigma$. Here 
$\sigma(s,t):p \rightarrow p \cdot s+t$ is a similarity of $\bR$, i.e. multiplication by $0 \ne s \in \bR$ and then addition by $t \in \bR$ for every 
$p \in \bR$. 
\end{definition}
\begin{rmrk}
If $\Gamma_1$ and $\Gamma_2$ are equivariant space groups then the their factor groups $\Gamma_1/L_{\Gamma_1}$ and $\Gamma_2/L_{\Gamma_2}$ are also
equivariant.  
\end{rmrk}
Thus the structure of the space group remains invariant under a similarity in the $\bR$-component and the spherical part is uniquely determined up to an isometry
of $\bS^2$. 

We characterize the spherical plane groups by the {\it Macbeath-signature} (see \cite{M}, \cite{Sz11}).

{\it In this paper we deal with such a $\SXR$ space group where the generators $\bg_i, \ (i=1,2,\dots m)$ of its point group $\Gamma_0$ 
are reflections and at least one of the possible translation parts of the above generators unequal to zero. 
These groups are called glide reflection groups.} \vspace{3mm}

\begin{rmrk} In \cite{Sz11} we have introduced the notion of {\it{generalized Coxeter group}}, if the generators $\bg_i, \ (i=1,2,\dots m)$ 
of its point group $\Gamma_0$ 
are reflections with translation parts $\tau_i=0, \ (i=1,2,\dots m).$
\end{rmrk}

In this paper we deal with the glide reflection space groups in $\SXR$ space which are by denotation of \cite{F01}:
\begin{enumerate}
\item $(+,~0,~[~~]~ \{(q,q)\}) \times \mathbf{1_R},~ q \ge 2$, \\ $\Gamma_0=(\bg_1,\bg_2 - \bg_1^2,\bg_2^2, (\bg_1\bg_2)^q)$,   

{\bf 2q.~I.~2}: $\big(\frac{1}{2},\frac{1}{2}\big)$;  {\bf 2qe.~I.~3}: $\big(0,\frac{1}{2}\big)$;
\item $(+,~0,~[~~]~ \{(2,2,q)\})\times \mathbf{1_R},~ q \ge 2$, \ \\ $\Gamma_0=(\bg_1,\bg_2,\bg_3 - \bg_1^2,\bg_2^2,\bg_3^2,(\bg_1\bg_3)^2, (\bg_2\bg_3)^2), 
(\bg_1\bg_2)^q) $,

{\bf 4q.~I.~2}: $\big(0,0,\frac{1}{2} \big)$;  {\bf 4q.~I.~3}: $\big(\frac{1}{2},\frac{1}{2},0 \big)$; 
{\bf 4q.~I.~4}: $\big(\frac{1}{2},\frac{1}{2}, \frac{1}{2} \big)$; \\ {\bf 4qe.~I.~5}: $\big(0, \frac{1}{2}, 0 \big)$;
{\bf 4qe.~I.~6}: $\big(0, \frac{1}{2}, \frac{1}{2} \big)$;
\item $(+,~0,~[~~]~ \{(2,3,3)\})\times \mathbf{1_R},$ \ \\ $\Gamma_0=(\bg_1,\bg_2,\bg_3 - \bg_1^2,\bg_2^2,\bg_3^2,(\bg_1\bg_2)^2, (\bg_1\bg_3)^3, (\bg_2\bg_3)^3), 
$ 

{\bf 11.~I.~2}: $\big(\frac{1}{2},\frac{1}{2}, \frac{1}{2} \big)$;
\item $(+,~0,~[~~]~ \{(2,3,4)\})\times \mathbf{1_R},$ \ \\ $\Gamma_0=(\bg_1,\bg_2,\bg_3 - \bg_1^2,\bg_2^2,\bg_3^2,(\bg_1\bg_2)^2, (\bg_1\bg_3)^3, (\bg_2\bg_3)^4), 
$ 

{\bf 12.~I.~2}: $\big(0, \frac{1}{2}, 0 \big)$;  {\bf 12.~I.~3}: $\big(\frac{1}{2}, 0, \frac{1}{2} \big)$; 
{\bf 12.~I.~4}: $\big(\frac{1}{2},\frac{1}{2}, \frac{1}{2} \big)$;
\item $(+,~0,~[~~]~ \{(2,3,5)\})\times \mathbf{1_R},$ \ \\ $\Gamma_0=(\bg_1,\bg_2,\bg_3 - \bg_1^2,\bg_2^2,\bg_3^2,(\bg_1\bg_2)^2, (\bg_1\bg_3)^3, (\bg_2\bg_3)^5), 
$ 

{\bf 13.~I.~2}: $\big(\frac{1}{2},\frac{1}{2}, \frac{1}{2} \big)$;
\end{enumerate}

\section{Geodesic curve and balls in $\SXR$ space}

E. {Moln\'ar} has shown in \cite{M97}, that the homogeneous 3-spaces
have a unified interpretation in the projective 3-sphere $\mathcal{PS}^3(\bV^4,\BV_4, \mathbb{R})$. 
In our work we shall use this projective model of $\SXR$ and
the Cartesian homogeneous coordinate simplex $E_0(\be_0)$,$E_1^{\infty}(\be_1)$,$E_2^{\infty}(\be_2)$,
$E_3^{\infty}(\be_3)$, $(\{\be_i\}\subset \bV^4$ \ $\text{with the unit point}$ $E(\be = \be_0 + \be_1 + \be_2 + \be_3 ))$ 
which is distinguished by an origin $E_0$ and by the ideal points of coordinate axes, respectively. 
Moreover, $\by=c\bx$ with $0<c\in \mathbb{R}$ (or $c\in\mathbb{R}\setminus\{0\})$
defines a point $(\bx)=(\by)$ of the projective 3-sphere $\cP \cS^3$ (or that of the projective space $\cP^3$ where opposite rays
$(\bx)$ and $(-\bx)$ are identified). 
The dual system $\{(\Be^i)\}\subset \BV_4$ describes the simplex planes, especially the plane at infinity 
$(\Be^0)=E_1^{\infty}E_2^{\infty}E_3^{\infty}$, and generally, $\Bv=\Bu\frac{1}{c}$ defines a plane $(\Bu)=(\Bv)$ of $\cP \cS^3$
(or that of $\cP^3$). Thus $0=\bx\Bu=\by\Bv$ defines the incidence of point $(\bx)=(\by)$ and plane
$(\Bu)=(\Bv)$, as $(\bx) \text{I} (\Bu)$ also denotes it. Thus {$\SXR$} can be visualized in the affine 3-space $\bA^3$
(so in $\bE^3$) as well.

In this context E. Moln\'ar \cite{M97} has derived the well-known infinitezimal arc-length square at any point of $\SXR$ as follows
\begin{equation}
   \begin{gathered}
     (ds)^2=\frac{(dx)^2+(dy)^2+(dz)^2}{x^2+y^2+z^2}.
       \end{gathered} \tag{2.1}
     \end{equation}
We shall apply the usual geographical coordiantes $(\phi, \theta), ~ (-\pi < \phi \le \pi, ~ -\frac{\pi}{2}\le \theta \le \frac{\pi}{2})$ 
of the sphere with the fibre coordinate $t \in \bR$. We describe points in the above coordinate system in our model by the following equations: 
\begin{equation}
x^0=1, \ \ x^1=e^t \cos{\phi} \cos{\theta},  \ \ x^2=e^t \sin{\phi} \cos{\theta},  \ \ x^3=e^t \sin{\theta} \tag{2.2}.
\end{equation}
Then we have $x=\frac{x^1}{x^0}=x^1$, $y=\frac{x^2}{x^0}=x^2$, $z=\frac{x^3}{x^0}=x^3$, i.e. the usual Cartesian coordinates.
We obtain by \cite{M97} that in this parametrization the infinitezimal arc-length square 
at any point of $\SXR$ is the following
\begin{equation}
   \begin{gathered}
      (ds)^2=(dt)^2+(d\phi)^2 \cos^2 \theta +(d\theta)^2.
       \end{gathered} \tag{2.3}
     \end{equation}
The geodesic curves of $\SXR$ are generally defined as having locally minimal arc length between their any two (near enough) points. 
The equation systems of the parametrized geodesic curves $\gamma(t(\tau),\phi(\tau),\theta(\tau))$ in our model can be determined by the 
general theory of Riemann geometry (see \cite{Sz11}).

Then by (2.2) we get with $c=\sin{v}$, $\omega=\cos{v}$ the equation systems of a geodesic curve, visualized in Fig.~1 in our Euclidean model:
\begin{equation}
  \begin{gathered}
   x(\tau)=e^{\tau \sin{v}} \cos{(\tau \cos{v})}, \\ 
   y(\tau)=e^{\tau \sin{v}} \sin{(\tau \cos{v})} \cos{u}, \\
   z(\tau)=e^{\tau \sin{v}} \sin{(\tau \cos{v})} \sin{u},\\
   -\pi < u \le \pi,\ \ -\frac{\pi}{2}\le v \le \frac{\pi}{2}. \tag{2.4}
  \end{gathered}
\end{equation}
\begin{rmrk}
Thus we have harmonized the scales along the fibre lines.
\end{rmrk}
\begin{definition}
The distance $d(P_1,P_2)$ between the points $P_1$ and $P_2$ is defined by the arc length of the shortest geodesic curve 
from $P_1$ to $P_2$.
\end{definition}
 \begin{definition}
 The geodesic sphere of radius $\rho$ (denoted by $S_{P_1}(\rho)$) with centre at the point $P_1$ is defined as the set of all points 
 $P_2$ in the space with the condition $d(P_1,P_2)=\rho$. Moreover, we require that the geodesic sphere is a simply connected 
 surface without selfintersection 
 in $\SXR$ space.
 \end{definition}
 \begin{rmrk}
 We shall see that this last condition depends on radius $\rho$.
 \end{rmrk}
 \begin{definition}
 The body of the geodesic sphere of centre $P_1$ and of radius $\rho$ in the $\SXR$ space is called geodesic ball, denoted by $B_{P_1}(\rho)$,
 i.e. $Q \in B_{P_1}(\rho)$ iff $0 \leq d(P_1,Q) \leq \rho$.
 \end{definition}
In \cite{Sz11} we have proved that $S(\rho)$ is a simply connected surface
in $\mathbf{E}^3$ if and only if $\rho \in [0,\pi)$, because if $\rho\ge\pi$ then there is at least one $v \in 
[-\frac{\pi}{2},\frac{\pi}{2}]$ so that $y(\tau,v)=z(\tau,v)=0$, i.e. selfintersection would occur (see (2.4)).
Thus we obtain the following
\begin{proposition}
The geodesic sphere and ball of radius $\rho$ exists in the $\SXR$ space if and only if $\rho \in [0,\pi].$
\end{proposition}
We have obtained (see \cite{Sz11}) the volume formula of the geodesic ball $B(\rho)$ of radius $\rho$ by the metric tensor $g_{ij}$ and by the 
Jacobian of (2.4):
\begin{Theorem}
\begin{equation}
\begin{gathered}
Vol(B(\rho))=\int_{V} \frac{1}{(x^2+y^2+z^2)^{3/2}}\mathrm{d}x ~ \mathrm{d}y ~ \mathrm{d}z = \\ = \int_{0}^{\rho} \int_{-\frac{\pi}{2}}^{\frac{\pi}{2}} 
\int_{-\pi}^{\pi} 
|\tau \cdot \sin(\cos(v)\tau)| ~ \mathrm{d} u \ \mathrm{d}v \ \mathrm{d}\tau = \\ =
2 \pi \int_{0}^{\rho} \int_{-\frac{\pi}{2}}^{\frac{\pi}{2}} |\tau \cdot \sin(\cos(v)\tau)| ~ \mathrm{d} v \ \mathrm{d}\tau. \tag{2.5}
\end{gathered}
\end{equation}
\end{Theorem}

\subsection{On fundamental domains}

A type of the fundamental domain of a studied space group can be combined as a fundamental domain of the corresponding spherical 
group with a part of a real 
line segment. This domain is called $\SXR$ {\it prism} (see \cite{Sz11}). 
{\it This notion will be important to compute the volume of the Dirichlet-Voronoi cell of a given space group because their volumes are equal and 
the volume of a $\SXR$ prism can be calculated by Theorem 2.8.}
\begin{figure}[ht]
\centering
\includegraphics[width=6cm]{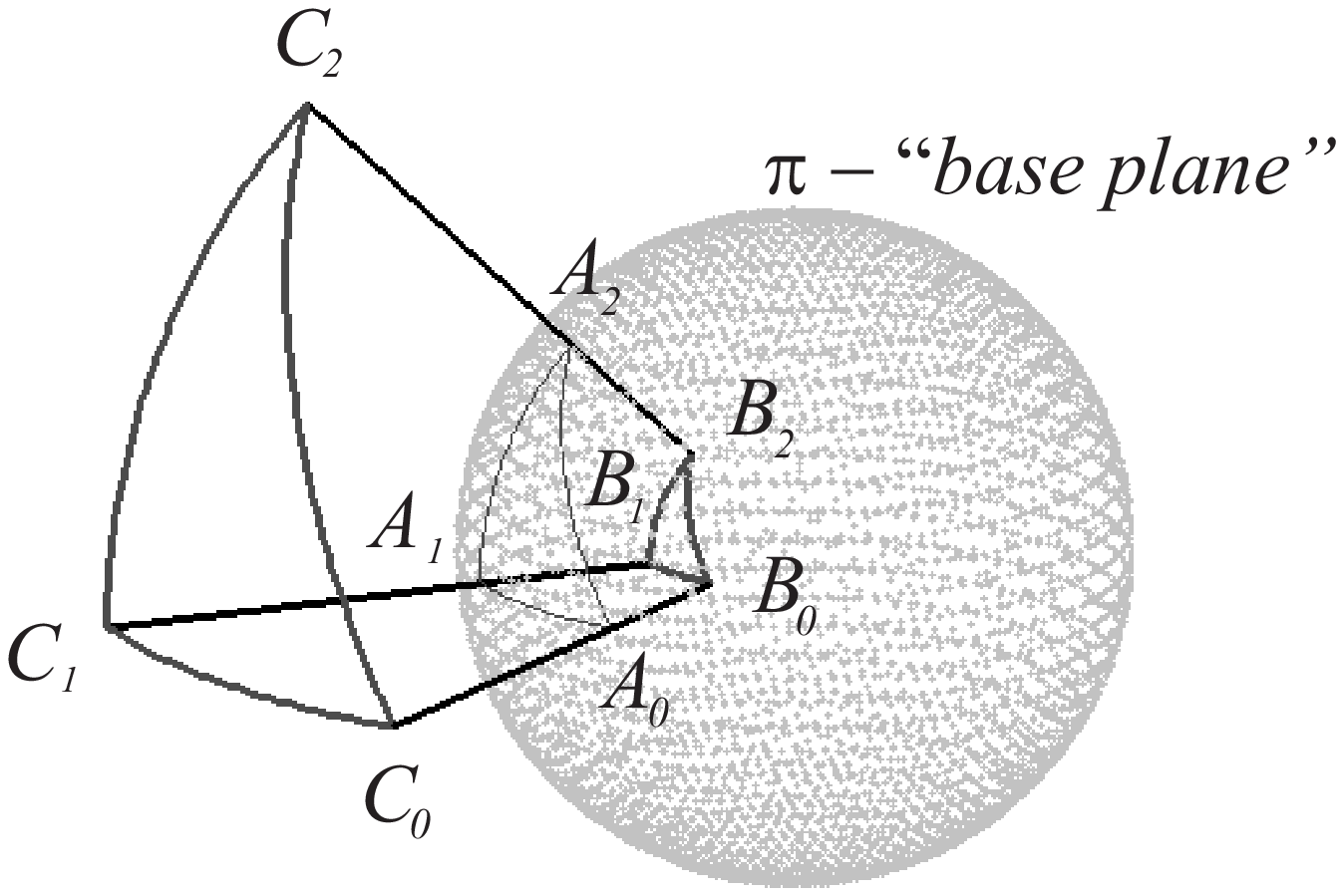} \includegraphics[width=4cm]{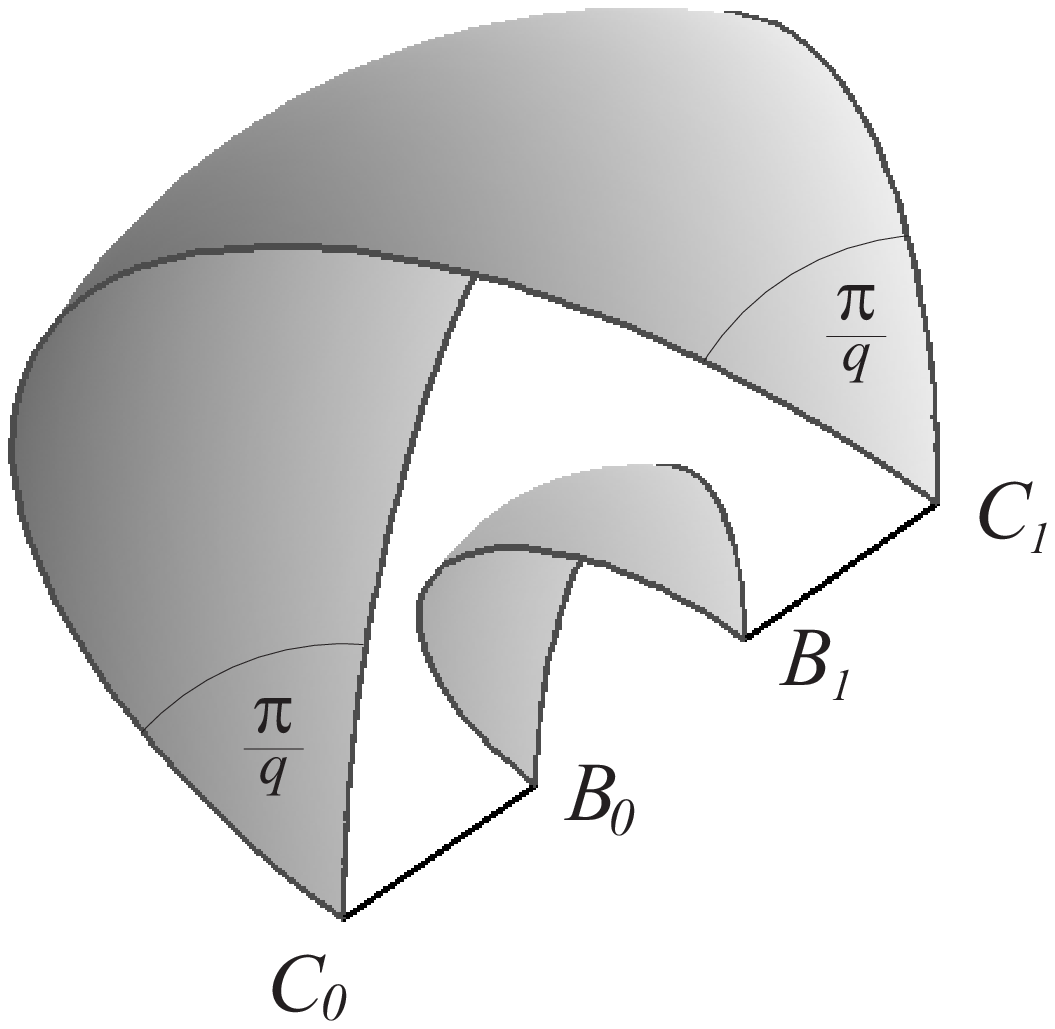}
\caption{Prism-like fundamental domains}
\label{}
\end{figure}

The $p$-gonal faces of a prism called cover-faces, and the other faces 
are the side-faces. The midpoints of the side edges 
form a "spherical plane" denoted by $\Pi$. It can be assumed that the plane $\Pi$ is the {\it base plane}:   
in our coordinate system (see (2.2)) the fibre coordinate $t=0$.

From \cite{Sz11} we recall
\begin{Theorem}
The volume of a $\SXR$ trigonal prism $\mathcal{P}_{B_0B_1B_2C_0C_1C_2}$ and of a digonal prism 
$\mathcal{P}_{B_0B_1C_0C_1}$ in $\SXR$ (see Fig.~1.a-b) can be computed by the following formula:
\begin{equation}
Vol(\mathcal{P})=\mathcal{A} \cdot h \tag{2.5}
\end{equation}
where $\mathcal{A}$ is the area of the spherical triangle $A_0A_1A_2$ or digon $A_0A_1$ in the base plane $\Pi$ with fibre 
coordinate $t=0$, and
$h=B_0C_0$ is the height of the prism.
\end{Theorem}
\section{Ball packings}
By remark (1.2) a $\SXR$ space group $\Gamma$ has a compact fundamental domain.
Usually the shape of the fundamental domain of a group of $\bS^2$ is not 
determined uniquely but the area of the domain is finite and unique 
by its combinatorial measure. Thus the shape of the fundamental domain of a crystallographic group of $\SXR$ is not unique as well.

In the following let $\Gamma$ be a fixed glide reflection space group of $\SXR$. We 
will denote by $d(X,Y)$ the distance of two points $X$, $Y$ by definition (2.2). 
\begin{definition}
We say that the point set
$$
\cD(K)=\{X\in\SXR\,:\,d(K,X)\leq d(K^\bg,X)\text{ for all }\bg\in\ \Gamma\}
$$
is the \emph{Dirichlet--Voronoi cell} (D-V~cell) to $\Gamma$ around the kernel 
point $K\in\SXR$.
\end{definition}
\begin{definition}
We say that
$$
\Gamma_X=\{\bg\in\Gamma\,:\,X^\bg=X\}
$$
is the \emph{stabilizer subgroup} of $X\in\SXR$ in $\Gamma$.
\end{definition}
\begin{definition}
Assume that the stabilizer $\Gamma_K=\bI$ i.e. $\Gamma$ acts simply transitively on 
the orbit of a point $K$. Then let $\cB_K$ denote the \emph{greatest ball} 
of centre $K$ inside the D-V cell $\cD(K)$, moreover let $\rho(K)$ denote the 
\emph{radius} of $\cB_K$. It is easy to see that
$$
\rho(K)=\min_{\bg\in\Gamma\setminus\bI}\frac12 d(K,K^\bg).
$$
\end{definition}
The $\Gamma$-images of $\cB_K$ form a ball packing $\cB^\Gamma_K$ with centre 
points $K^\bG$. 
\begin{definition}
The \emph{density} of ball packing $\cB^\Gamma_K$ is
$$
\delta(K)=\frac{Vol(\cB_K)}{Vol\cD(K)}.
$$
\end{definition}
It is clear that the orbit $K^\Gamma$ and the ball packing $\cB^\Gamma_K$ have the 
same symmetry group, moreover this group contains the starting 
crystallographic group $\Gamma$:
$$
Sym K^\Gamma=Sym\cB^\Gamma_K\geq\Gamma.
$$
\begin{definition}\rm
We say that the orbit $K^\Gamma$ and the ball packing $\cB^\Gamma_K$ is 
\emph{characteristic} if $Sym K^\Gamma=\Gamma$, else the orbit is not 
characteristic.
\end{definition}

\subsection{Simply transitive ball packings}

\emph{Our problem is} to find a 
point $K\in\ \SXR$ and the orbit $K^\Gamma$ for $\Gamma$ such that $\Gamma_K=\bI$ 
and the density $\delta(K)$ of the corresponding ball packing 
$\cB^\Gamma(K)$ is maximal. In this case the ball packing $\cB^\Gamma(K)$ is 
said to be \emph{optimal.} 

The lattice of $\Gamma$ has a free parameter $p(\Gamma)$. Then we have to find the densest ball packing on $K$ for fixed 
$p(\Gamma)$, and vary $p$ to get the optimal ball packing. 
\begin{equation}
\delta(\Gamma)=\max_{K, \ p(\Gamma)}(\delta(K)) \tag{3.1}
\end{equation}
Let $\Gamma$ be a fixed {\it glide reflection group}. 
The stabiliser of $K$ is trivial i.e. we are looking the optimal kernel point
in a 3-dimensional region, inside of a fundamental domain of $\Gamma$ with free fibre parameter $p(\Gamma)$.
It can be assumed by the homogeneity of $\SXR$, that the fibre coordinate of the center of the optimal ball is zero.

\subsection{Optimal ball packing to space group {\bf 12.~I.~3}}

Now we consider the following point group:

\begin{equation}
\begin{gathered}
(+,0;[~]; ~ \{(2,3,4)\}) \times \mathbf{1}_\bR; \\ \Gamma_0:=\{\bg_1,\bg_2,\bg_3 ~-~ \bg_1^2,\bg_2^2,\bg_3^2, (\bg_1 \bg_2)^2, (\bg_1 \bg_3)^3, (\bg_2 \bg_3)^4\}.   
\end{gathered} \tag{3.2}
\end{equation}
This is the full isometry group of the usual cube surface, generated by the three reflections $\bg_i, ~ i=1,2,3$. 
The possible translation parts of the generators of $\Gamma_0$ will be determined by (1.2) and by the defining relations of the point group. 
Finally, 
from the so-called Frobenius congruence relations we obtain the four non equivariant solutions:
$$(\tau_1,\tau_2,\tau_3) \cong (0,0,0),~(0,\frac{1}{2},0),~(\frac{1}{2},0,\frac{1}{2}),~(\frac{1}{2},\frac{1}{2},\frac{1}{2}).$$
If $(\tau_1,\tau_2,\tau_3) \cong ~(\frac{1}{2},0,\frac{1}{2})$ then we get the $\SXR$ space group {\bf 12.~I.~3}.
The fundamental domain of its point group is a spherical triangle $A_0A_1A_2$ with 
angles $\frac{\pi}{3}$, $\frac{\pi}{2}$, $\frac{\pi}{4}$ lying in the base plane $\Pi$ (see Fig.~2).
It can be assumed by the homogeneity of $\SXR$, that the fibre coordinate of the center of the optimal ball is zero and it is an 
interior point of $A_0A_1A_2$ triangle.

We shall apply the Cartesian homogeneous coordinate system introduced in Section 2 (see Fig.~2) and the usual geographical coordiantes 
$(\phi, \theta), ~ (-\pi < \phi \le \pi, ~ -\frac{\pi}{2}\le \theta \le \frac{\pi}{2})$ 
of the sphere with the fibre coordinate $t \in \bR$ (see (2.2)).

We consider an arbitrary interior point $K(x^0,x^1,x^2,x^3)$ of spherical triangle $A_0A_1A_2$  
in the above coordinate system in our model by the following equations: 
\begin{equation}
x^0=1, \ \ x^1= \cos{\phi} \cos{\theta},  \ \ x^2=\sin{\phi} \cos{\theta},  \ \ x^3= \sin{\theta} \tag{3.3}
\end{equation}
\begin{figure}[ht]
\centering
\includegraphics[width=9cm]{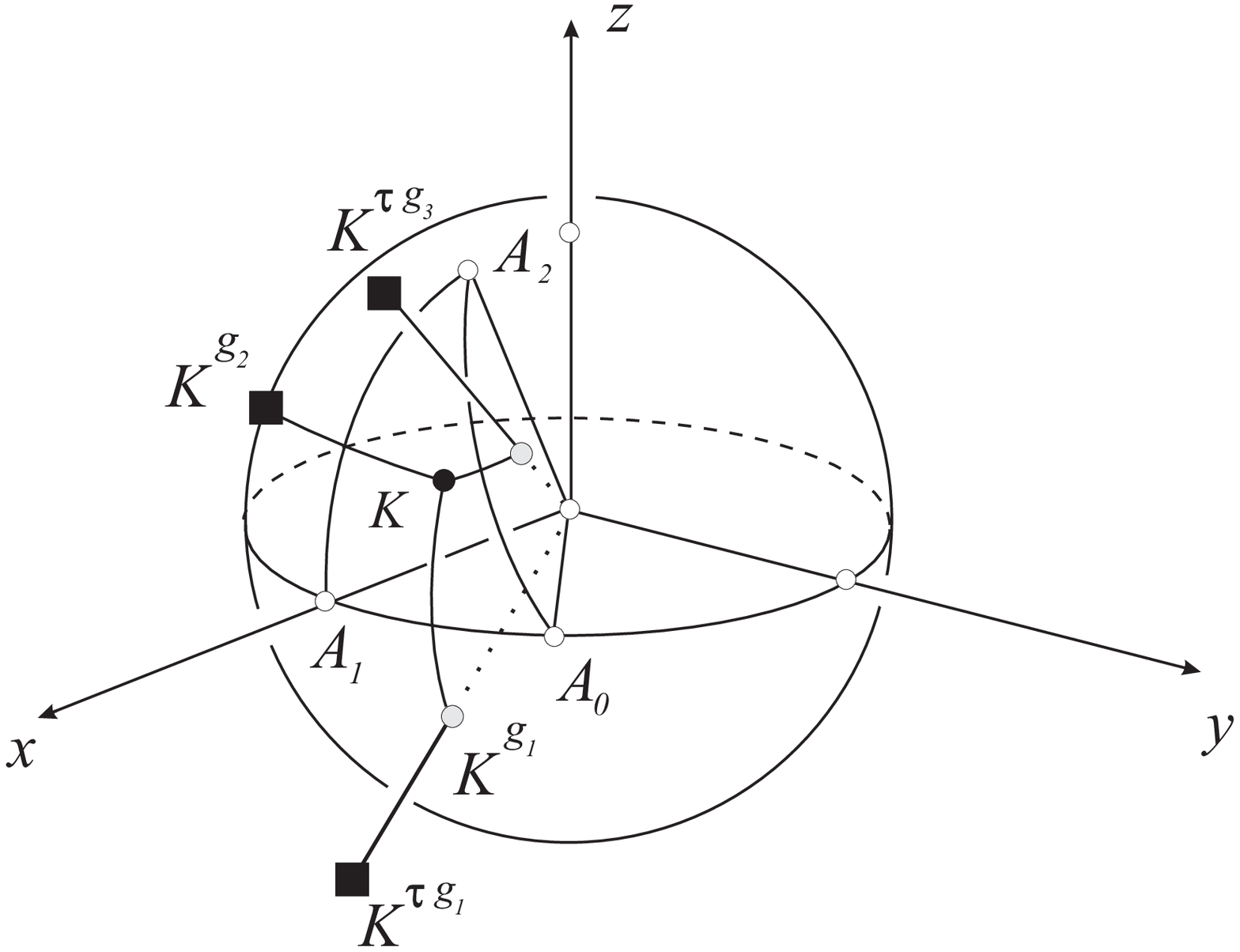} 
\caption{}
\label{}
\end{figure}
Let $\mathcal{B}_\Gamma(R)$ denote a geodesic ball packing of $\SXR$ space with balls $B(R)$ of radius $R$ where their 
centres give rise to the orbit $K^\Gamma$. In the following we consider to each ball packing the possible smallest 
translation part $\tau(K,R)$ (see Fig.~2) depending on $\Gamma$, $K$ and $R$. 
A fundamental domain of $\Gamma$ is its D-V cell $\cD(K)$ around the kernel point $K$. It is clear that the optimal ball $\mathcal{B}_K$ has to touch 
some faces of its D-V cell. The volume of $D(K)$ is equal to the volume of the prism which is given 
by the fundamental domain of the point group $\Gamma_0$ of $\Gamma$ and by the height $2\tau(R,K)$. The images of $D(K)$ 
by our discrete isometry group covers the $\SXR$ space without overlap. 
For the density of the packing it is sufficient to relate the volume of the optimal ball
to that of the solid $D(K)$ (see Definition 3.4).

It is clear, that the densest ball arrangement $\mathcal{B}_\Gamma(R)$ of balls $B(R)$ has to hold the following requirements:
\begin{equation}
\begin{gathered}
(a) \ \ d(K,K^{g_2})=2R=d(K,K^{\tau g_1}), \\
(b) \ \ d(K,K^{g_2})=2R=d(K,K^{\tau g_3}), \\
(c) \ \ d(K,K^{2\tau}) \ge 2R \\
(d) \ \ \text{Balls of radius $R$ with centres} \\ \text{$K$, $K^{g_2}$, $K^{\tau g_1}$, $K^{\tau g_3}$, $K^{2\tau}$ form a packing.}
\end{gathered} \tag{3.5}
\end{equation}

Here $d$ is the distance function in the $\SXR$ space (see Definition 2.2). The equations (a) and (b) mean that the ball centres $K^{\tau g_1}$ and 
$K^{\tau g_3}$ 
lie on the equidistant geodesic surface of the points $K$ and $K^{2\tau}$ which is a sphere in our model in this case (see \cite{PSSz}). 

We consider two main ball arrengements: 
\begin{enumerate}
\item We denote by $\mathcal{B}_{\Gamma}(R_0,K_0)$ those packing where requirements (3.5)  
and \\ $d(K,K^{2\tau}) = 2R$ hold (see Fig.~3).
\item We denote by $\mathcal{B}_{\Gamma}(R_1,K_1)$ those packing where requirements (3.5)  
and  \\ $d(K^{\tau g_1},K^{\tau g_3}) =2R$ hold (see Fig.~4).
\end{enumerate}
\begin{figure}
\centering
\includegraphics[width=6cm]{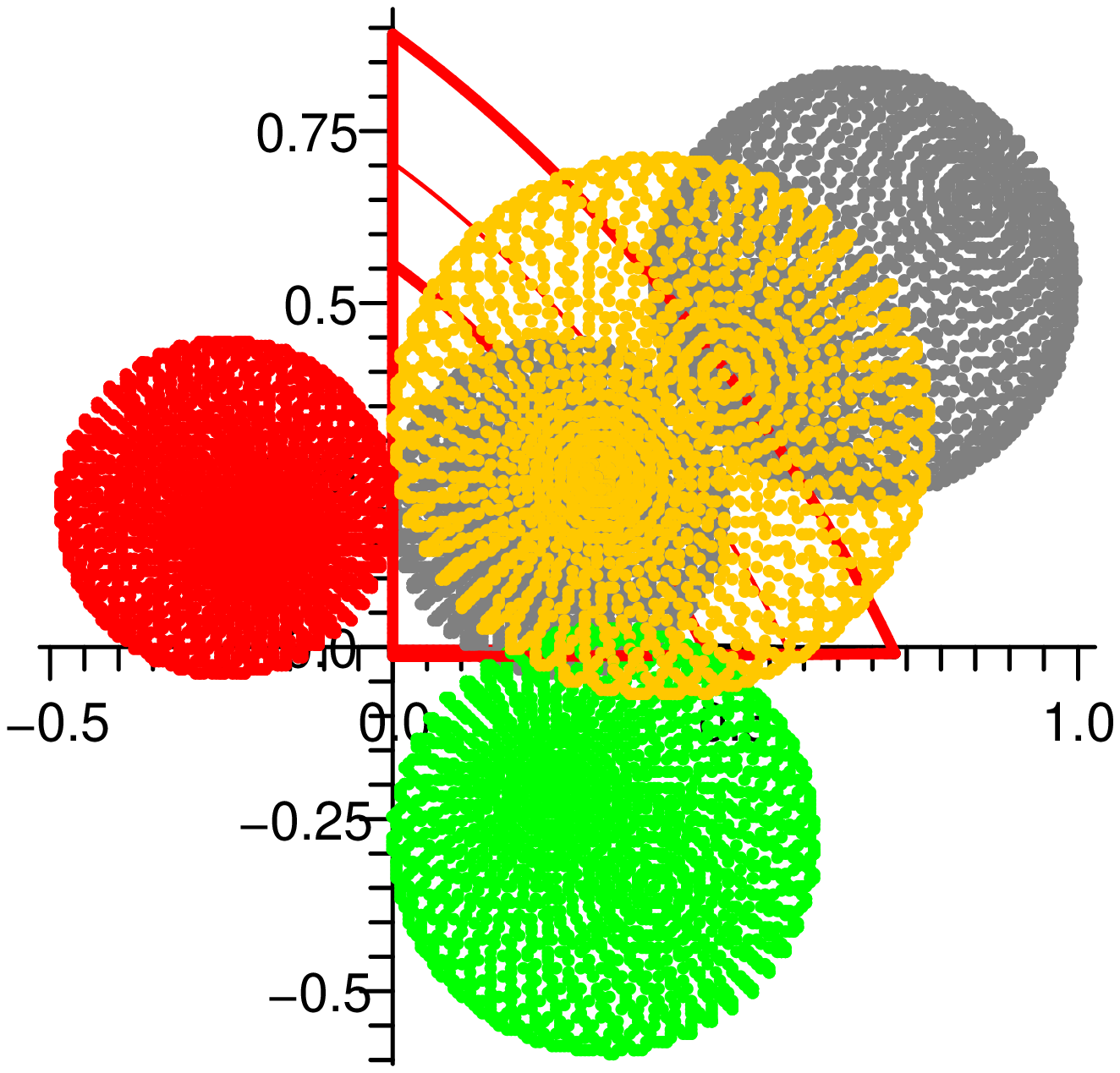} \includegraphics[width=6cm]{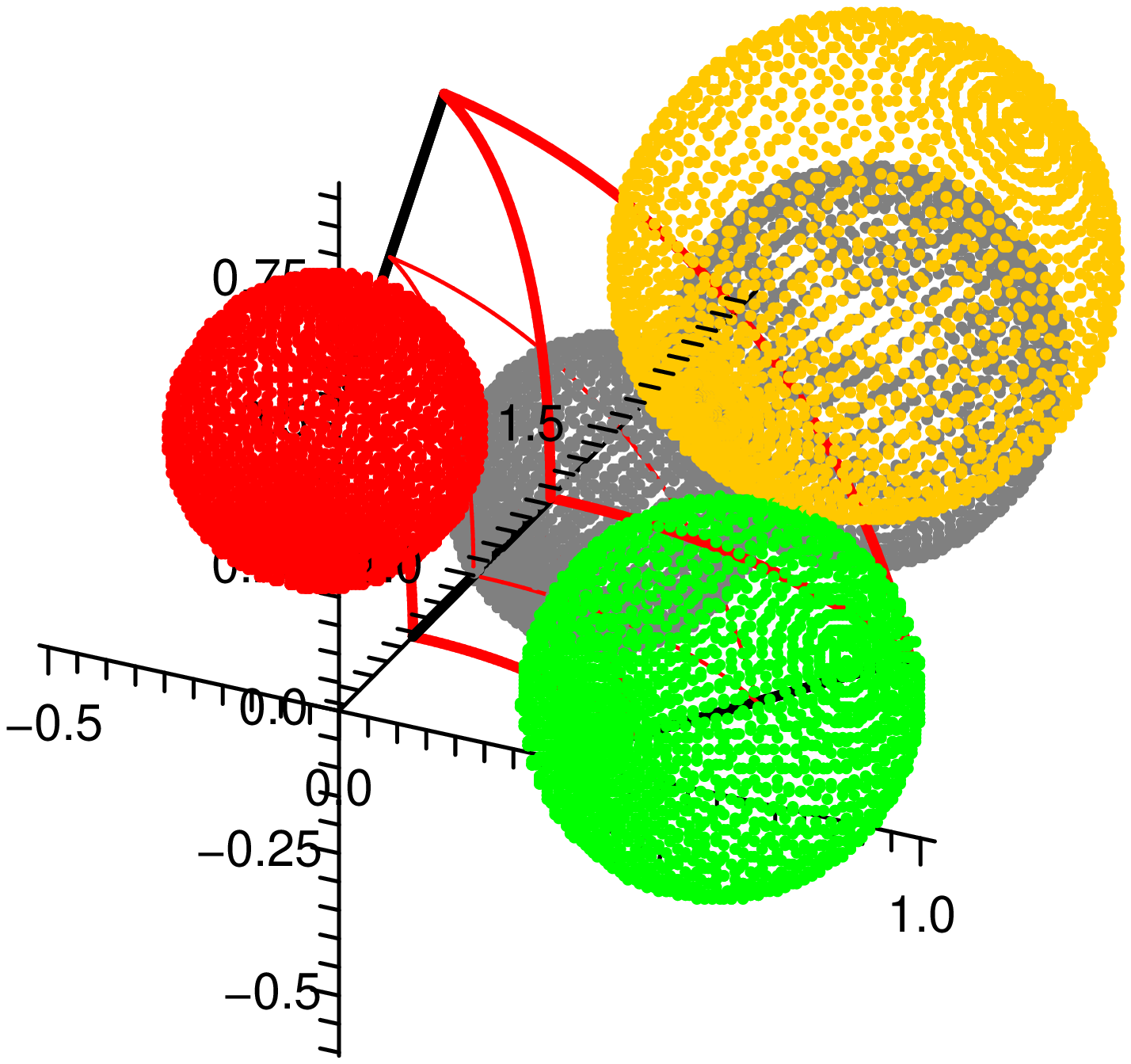}
\caption{}
\label{}
\end{figure}
First we determine the coordinates of the points $K_i$, ($i=1,~2$) ($K_i$ is given by (3.3) with parameters $\phi$ and $\theta$), the radius $R$ of the ball, 
the volume of a ball $B(R)$  and the density of the packing in both main cases. 
We get the following solutions by systematic approximation, where the computations were carried out by {\it Maple V Release 10} up to 30 decimals:
\begin{equation}
\begin{gathered}
\phi_0 \approx 0.24389626,\ \ \theta_0 \approx 0.20663860, \ \ R_0 \approx 0.23860571, \\
Vol(B(R_0)\approx 0.05668684, \ \ \delta(R_0,K_0) \approx 0.45373556.
\end{gathered} \tag{3.5}
\end{equation}
\begin{equation}
\begin{gathered}
\phi_1 \approx 0.30773985,\ \ \theta_1 \approx 0.17313169, \ \ R_1 \approx 0.30299179, \\
Vol(B(R_1)\approx 0.11580359, \ \ \delta(R_1,K_1) \approx 0.44472930.
\end{gathered} \tag{3.6}
\end{equation}
We obtain by careful investigation of the density function 
$\delta(R,K)$ ($R \in [R_0,R_1]$) of the considered ball packing  the following:
\begin{Theorem}
The ball arrangement $\mathcal{B}_{\Gamma}(R_0,K_0)$ (see Fig.~3) provides the densest symply transitive ball packing belonging to the $\SXR$ 
space group {\bf 12.~I.~3}. 
\end{Theorem}
\begin{figure}
\centering
\includegraphics[width=6cm]{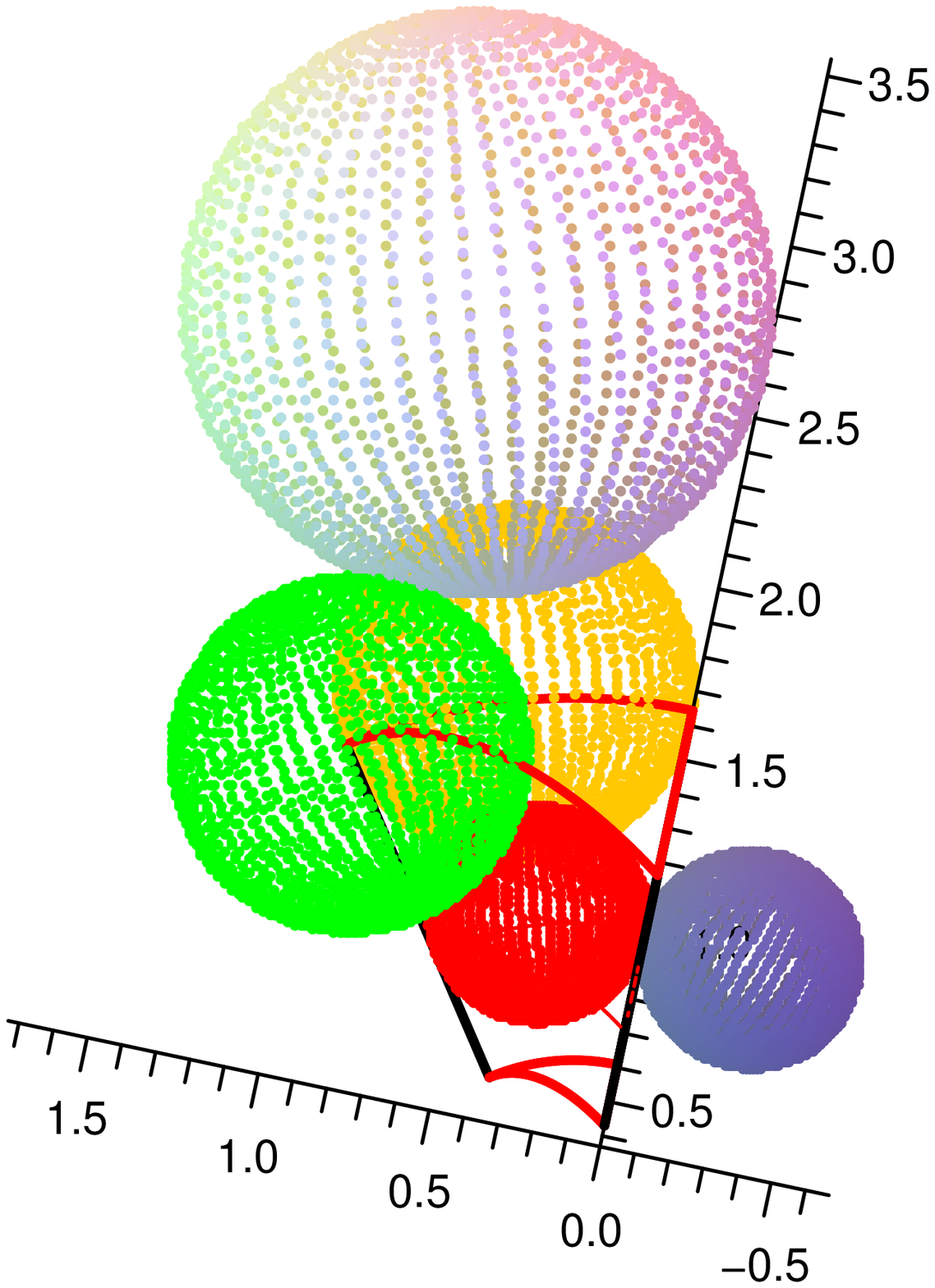} \includegraphics[width=6cm]{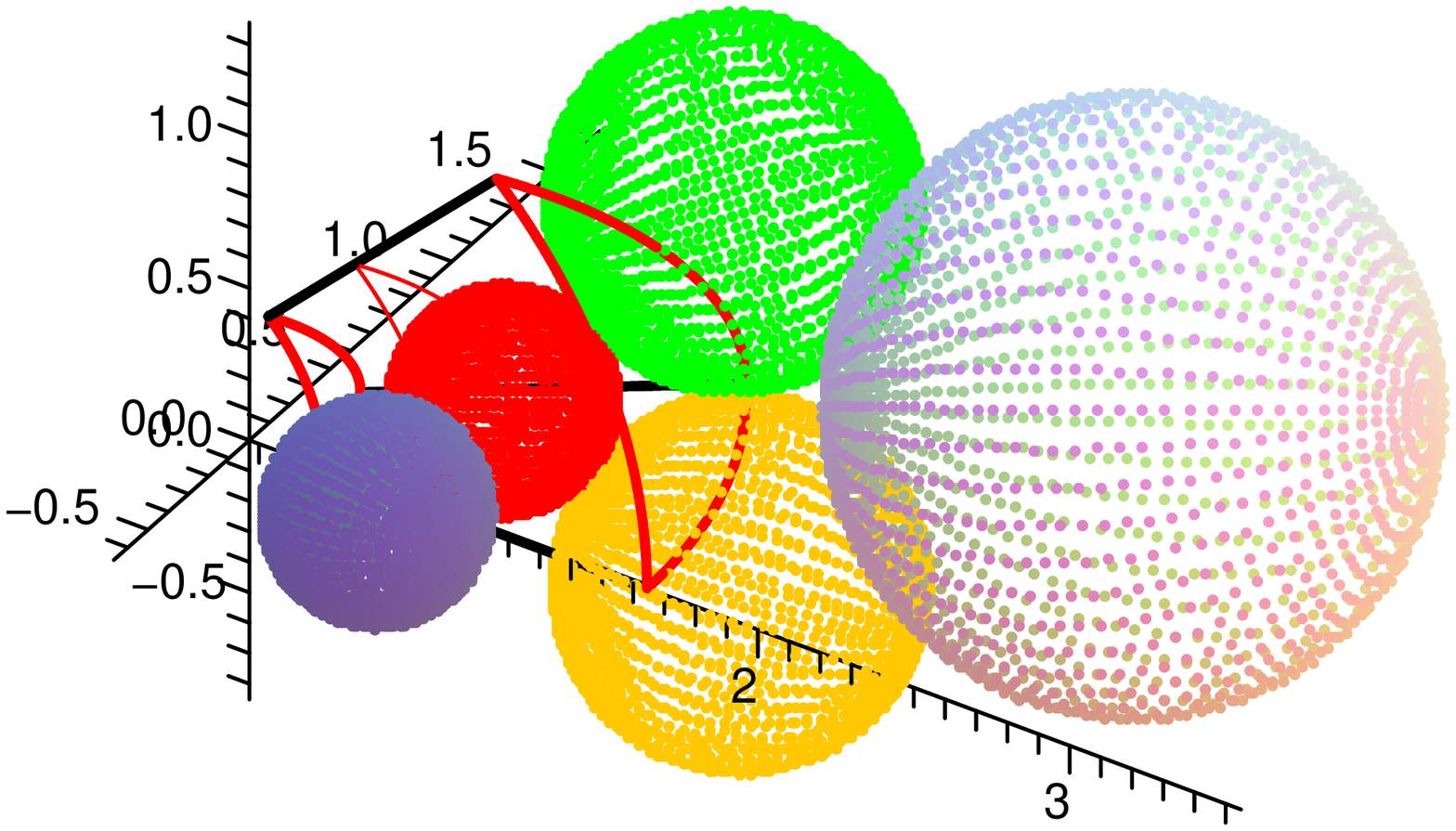}
\caption{}
\label{}
\end{figure}
\subsection{The densest simply transitive ball packing}
We consider the following point group:
\begin{equation}
\begin{gathered}
(+,~0,~[~~]~ \{(q,q)\}) \times \mathbf{1_R},~ q \ge 2; \\ \Gamma_0=(\bg_1,\bg_2 - \bg_1^2,\bg_2^2, (\bg_1\bg_2)^q).   
\end{gathered} \notag
\end{equation}
This point group is generated by two reflections $\bg_i, ~ i=1,2,3$. 
The possible translation parts of the generators of $\Gamma_0$ will be determined by (1.2) and by the defining relations of the point group. 
Finally, we obtain
from the so-called Frobenius congruence relations three non equivariant solutions:
$$(\tau_1,\tau_2) \cong (0,0),~(0,\frac{1}{2}),~(\frac{1}{2},\frac{1}{2}).$$
If $(\tau_1,\tau_2) \cong ~(\frac{1}{2},\frac{1}{2})$ then we have obtained the $\SXR$ space group {\bf 2q.~I.~2}.

The fundamental domain of the point group of the considered space group is a spherical digon $A_0A_1$ with 
angle $\frac{\pi}{q}$ in the base plane $\Pi$.
Similarly to the above section can be assumed, that the fibre coordinate of the center of the optimal ball is zero and it is an 
interior point of $A_0A_1$ digon (see Fig.~5).  

{\it In the following we consider ball packings belonging to $q=2$.}
\begin{figure}
\centering
\includegraphics[width=9cm]{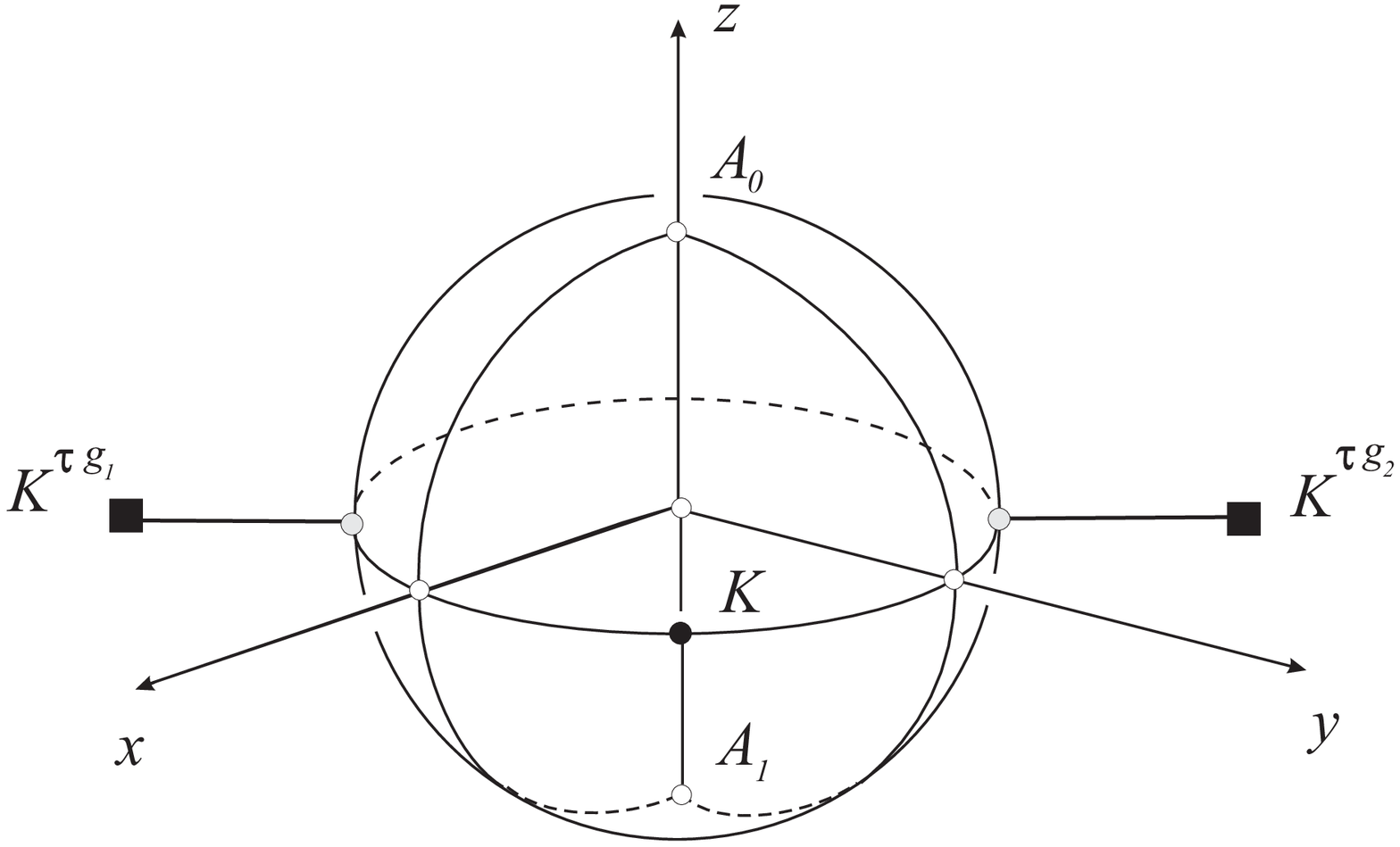}
\caption{}
\label{}
\end{figure}
We use also the above introduced Cartesian homogeneous coordinate system and the usual geographical coordinates 
$(\phi, \theta), ~ (-\pi < \phi \le \pi, ~ -\frac{\pi}{2}\le \theta \le \frac{\pi}{2})$ 
of the sphere with the fibre coordinate $t \in \bR$ (see (2.2)).

We consider an arbitrary interior point $K(1,x^1,x^2,x^3)=K(\phi,\theta)$ of spherical digon $A_0A_1$  
in the above coordinate system in our model (see Fig.~5).

Our aim is to determine the maximal radius $R$ of the balls, and the maximal density $\delta(K,R)$.

The ball arrangement $\mathcal{B}_{opt}(K,R)$ is defined by the following equations:
\begin{equation}
\begin{gathered}
(a) \ \ d(K,K^{\tau g_1})=2R=d(K,K^{\tau g_2}), \\
(b) \ \ d(K^{\tau g_1},K^{\tau g_2})=2R, \\
\end{gathered} \tag{3.7}
\end{equation}
We can determine the coordinates of the point $K$, the radius $R$ of the ball, 
the volume of a ball $B(R)$  and the density of this packing: 
\begin{figure}[ht]
\centering
\includegraphics[width=6cm]{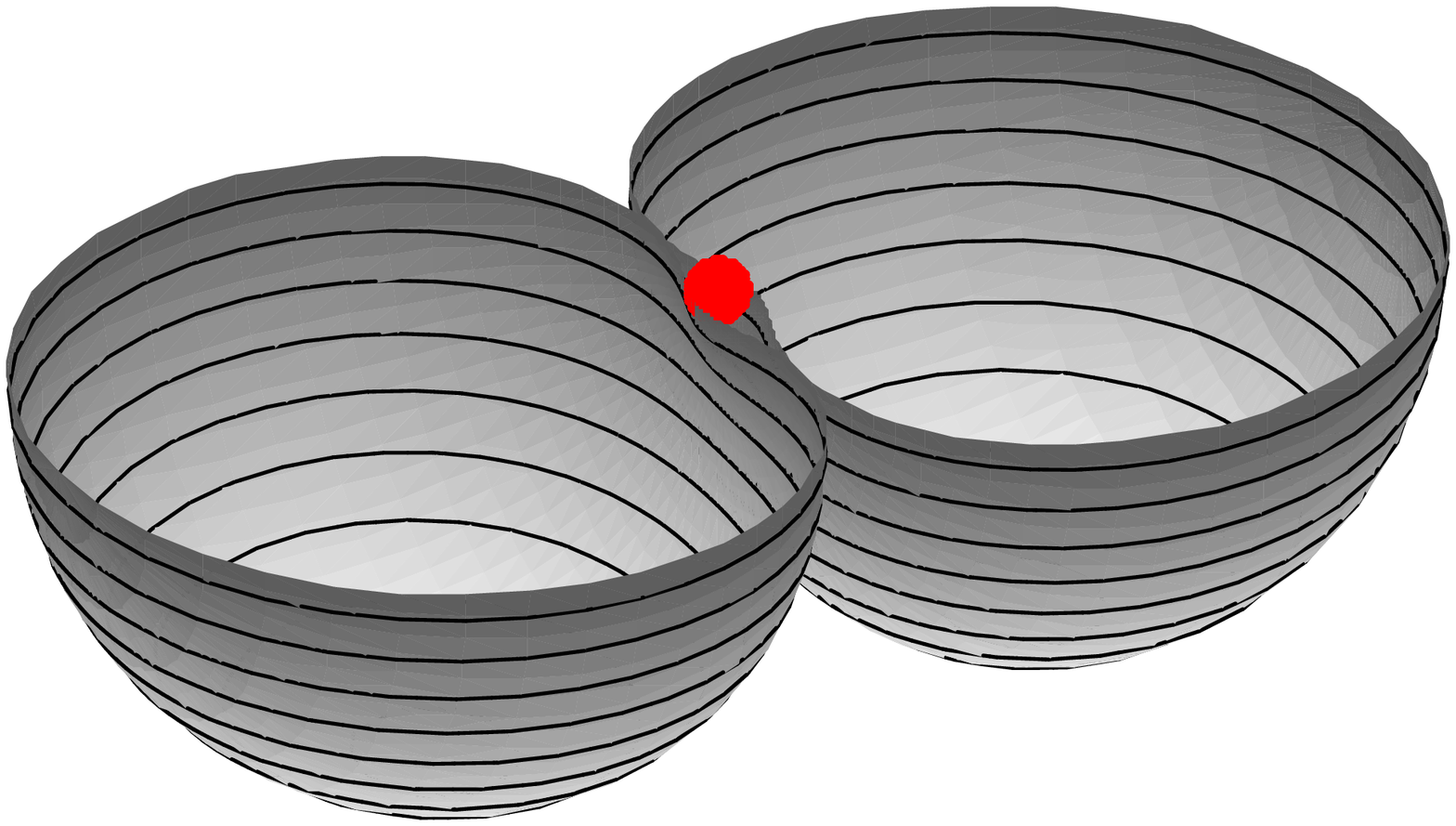} \includegraphics[width=6cm]{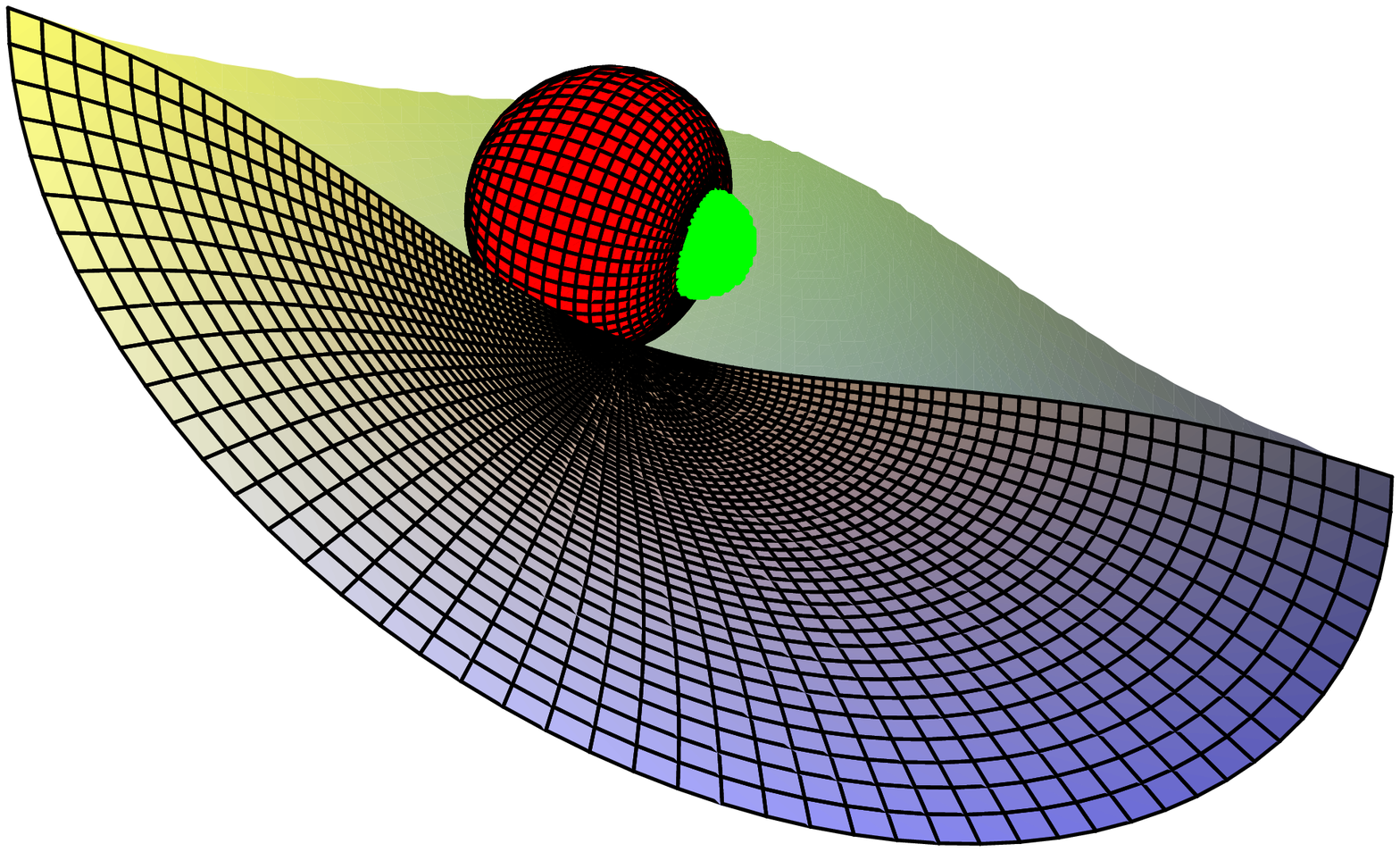}
\caption{The densest simply transitive geodesic ball packing}
\label{}
\end{figure}
\begin{equation}
\begin{gathered}
\phi= \frac{\pi}{4} \approx 0.78539816,\ \ \theta =0, \ \ R = \frac{\pi}{2}\approx 1.57079633, \\
Vol(B(R)\approx 13.74539472, \ \ \delta(R,K) \approx 0.80407553.
\end{gathered} \tag{3.8}
\end{equation}
Similarly to the above Section we can prove the following theorem:
\begin{Theorem}
The ball arrangement $\mathcal{B}_{opt}(R,K)$ (see Fig.~6) provides the densest simply transitive ball packing belonging to 
the $\SXR$ space group {\bf 2q.~I.~2}.
\end{Theorem}
Considering all in this paper investigated space groups and in \cite{Sz11} discussed generalized Coxeter $\SXR$ space groups we get the next
theorem:
\begin{Theorem}
The ball arrangement $\mathcal{B}_{opt}(R,K)$ provides the densest simply transitive ball packing belonging to 
the generalized Coxeter and glide reflections generated $\SXR$ space groups.
\end{Theorem}
By Theorems 2.7 and 2.8 and by Definitions 3.3 and 3.4, similarly to the above space groups, we have determined
the data (radii, densities and volumes of optimal balls) of the optimal simply transitive ball packings to each 
glide reflections generated $\SXR$ space group which are summarized in Table 1.
\begin{rmrk}
The space groups {\bf 2q.~I.~2}, {\bf 2qe.~I.~3}, {\bf 4q.~I.~2}, {\bf 4q.~I.~3}, {\bf 4q.~I.~4}, {\bf 4qe.~I.~5}, {\bf 4qe.~I.~6}
depend on parameter $q$ thus their optimal ball packings depend also on $q$ but in the Table 1 we give only the data of the densest 
ball packing indicating its $q$ parameter to each considered space group.   
\end{rmrk}
\medbreak
\centerline{\vbox{
\halign{\strut\vrule ~ \hfil $#$ \hfil  \vrule
& ~ \hfil $#$ \hfil \vrule
& ~ \hfil $#$ \hfil \vrule
& ~ \hfil $#$ \hfil \vrule
\cr
\noalign{\hrule}
\multispan4{\strut\vrule\hfill{\it  Table 1} \hfill\vrule}%
\cr
\noalign{\hrule}
\noalign{\vskip1pt}
\noalign{\hrule}
Space ~ group & R  & Vol(\mathcal{B}_{{K}}(R)) & \delta \cr
\noalign{\hrule}
\noalign{\vskip1pt}
\noalign{\hrule}
\noalign{\vskip1pt}
{\bf 2q.~I.~2}, ~ q=2 &  \frac{\pi}{2}\approx 1.57079633 & \approx 13.74539472 & \approx 0.80407553\cr
\noalign{\vskip1pt}
\noalign{\hrule}
\noalign{\vskip1pt}
{\bf 2qe.~I.~3}, ~ q=2 &  \frac{\pi}{2}\approx 1.57079633 & \approx 13.74539472 & \approx 0.69634983  \cr
\noalign{\vskip1pt}
\noalign{\hrule}
\noalign{\vskip1pt}
{\bf 4q.~I.~2},~ q=2 & \approx 0.64360446 & \approx 1.08624788 & \approx 0.53722971 \cr
\noalign{\vskip1pt}
\noalign{\hrule}
{\bf 4q.~I.~3}, ~ q=2& \approx 0.67517586 & \approx 1.25058159 & \approx 0.58958340 \cr
\noalign{\vskip1pt}
\noalign{\hrule}
\noalign{\vskip1pt}
{\bf 4q.~I.~4}, ~ q=2 & \approx 0.95531662 & \approx 3.43551438 & \approx 0.74837055 \cr
\noalign{\vskip1pt}
\noalign{\hrule}
\noalign{\vskip1pt}
{\bf 4qe.~I.~5}, ~ q=2 &  \approx 0.64360446 & \approx 1.08624788 & \approx 0.53722971 \cr
\noalign{\vskip1pt}
\noalign{\hrule}
\noalign{\vskip1pt}
{\bf 4qe.~I.~6}, ~ q=2& \approx 0.67517586 & \approx 1.25058159 & \approx 0.58958340 \cr
\noalign{\vskip1pt}
\noalign{\hrule}
{\bf 11.~I.~2}& \approx 0.46364761 & \approx 0.41154972 & \approx 0.58861600 \cr
\noalign{\vskip1pt}
\noalign{\hrule}
{\bf 12.~I.~2}& \approx 0.22770028 & \approx 0.04928081 & \approx 0.41334779 \cr
\noalign{\vskip1pt}
\noalign{\hrule}
{\bf 12.~I.~3}& \approx 0.23860571 & \approx 0.05668684 & \approx 0.45373556 \cr
\noalign{\vskip1pt}
\noalign{\hrule}
{\bf 12.~I.~4}&  \approx 0.31004511 & \approx 0.12404486 & \approx 0.53597559 \cr
\noalign{\vskip1pt}
\noalign{\hrule}
\noalign{\vskip1pt}
{\bf 13.~I.~2} &  \approx 0.18705243 & \approx 0.02735051 & \approx 0.49222087  \cr
\noalign{\vskip1pt}
\noalign{\hrule}
\noalign{\vskip1pt}
}}}
It is timely to arising the above question for further space groups in $\SXR$ space.

In this paper we have mentioned only some problems in discrete geometry of $\SXR$ space, but we hope that from these 
it can be seen that our projective
method suits to study and solve similar problems (\cite{PSSz}, \cite{Sz10}, \cite{Sz11}). Analogous questions in other homogeneous Thurston geometries are interesting.

\end{document}